\documentclass[12pt]{article}
\usepackage{amsfonts}
\usepackage{amssymb}
\usepackage{graphics}
\usepackage[dvips]{epsfig}

\def\cC{{\cal C}}

\def\cL{{\cal L}}

\def\d{\delta}
\def\e{\epsilon}
\def\z{\zeta}

\def\l{\lambda}

\def\s{\sigma}

\def\R{\mathbb{R}}
\def\C{\mathbb{C}}

\outer\def\proclaim#1:#2\par{\medbreak\vskip-\parskip
    \noindent{\bf#1.\enspace}{\sl#2}
  \ifdim\lastskip<\medskipamount \removelastskip\penalty55\medskip\fi}
\def\demo#1:{\par\medskip\noindent\it{#1}. \rm}
\def\endpr{\hfill $\spadesuit$ \medskip}
\def\ni{\noindent}               
\def\bar{\overline}              
\def\bs{\backslash}              
\def\di{\partial}                
\def\wt{\widetilde}
\def\holo{holomorphic}                   
\def\nbd{neighborhood}                   
\def\psc{pseudoconvex}                   
\def\spsc{strongly\ pseudoconvex}        
\def\ra{real-analytic}                   
\def\spsh{strongly\ plurisubharmonic}
\def\tr{totally real}                    
\def\iff{if and only if}

\title{STRONGLY PSEUDOCONVEX HANDLEBODIES}

\author{Franc Forstneri\v c and Jernej Kozak
        \date{}
       }

\begin{document}
\baselineskip=16pt
\thispagestyle{empty}

\maketitle
\begin{center} {\bf Contents} \end{center}
\bigskip
\begin{flushleft}
  1. Introduction \dotfill $\quad \pageref{introduction}$\\
  2. Pseudoconvexity of spherical domains               
     \dotfill $\quad \pageref{Pseudoconvexity}$\\
  3. Strongly pseudoconvex handlebodies             
     \dotfill $\quad \pageref{SPH}$\\
  4. Handlebodies on general quadratic domains                      
     \dotfill $\quad \pageref{HandlebodiesGQD}$\\
  References \dotfill $\,\,\, \pageref{References}$
\end{flushleft}
\vskip 8mm
\noindent
\section{Introduction}\label{introduction}
Let $\C^n$ be the complex $n$-dimensional Euclidean space with coordinates
$z=(z_1,\ldots, z_n)$, $z_j=x_j+i y_j$. Let $J$ denote the standard almost complex 
structure operator on $T\C^n$: $J({\di\over \di x_j})={\di \over \di y_j}$,
$J({\di\over \di y_j})=-{\di \over \di x_j}$. A $\cC^1$ submanifold 
$M\subset \C^n$ is {\it totally real\/} at $p\in M$ if $T_p M\cap JT_p M=\{0\}$, 
that is, the tangent space $T_p M \subset T_p\C^n$ contains no complex line. 
A $\cC^2$ function $\rho\colon U\subset\C^n\to \R$ is {\it \spsh\/}  on $U$ if 
$$ 
	\cL_\rho(z;v) = \sum_{j,k=1}^n 
	{\di^2\rho \over \di z_j\di\bar z_k}(z) \, v_j\bar v_k >0
        \quad (z\in U,\ v\in \C^n\bs \{0\}).                               
$$ 
$\cL_\rho(z;v)$ is called the {\it Levi form} of $\rho$ at $z$ in the direction of 
the vector $v$.

Assume that $D\subset \C^n$ is a closed, smoothly bounded, 
\spsc\ domain. Thus $D=\{\rho\le 0\}$ where $\rho$ is a \spsh\ function 
in an open set $U\supset D$, with $d\rho\ne 0$ on $bD=\{\rho=0\}$.
Let $M\subset \C^n$ be a smooth \tr\ submanifold with boundary
$bM=S\cup S'$, where each of the sets $S$, $S'$ is a union 
of connected components of $bM$ ($S'$ may be empty). Assume 
furthermore that 
$$
	M\cap D = S \subset bD,  \qquad 
	T_p (S) \subset T_p^\C (bD) := T_p(bD) \cap JT_p(bD) \quad (p\in S).   
$$
Such $M$ will be called a {\it totally real handle} attached to $D$ 
along the {\it Legendrian} ({\it complex tangential}) submanifold 
$S\subset bD$. (Some authors reserve the word `handle' for the case when $M$ 
is diffeomorphic to the closed ball in  some $\R^k$ and $bM=S^{k-1}$.) 
We consider the following problem.

%
%
\proclaim 
The handlebody problem: Given a (small) open set $U\supset M$, find a closed 
\spsc\ domain $K\subset\C^n$ satisfying $D\cup M\subset K  \subset D\cup U$ 
which admits a strong deformation retraction onto $E:=D\cup M$.

%
%
\begin{center}
\begin{figure}[hbt]
  \centering\epsfig{file=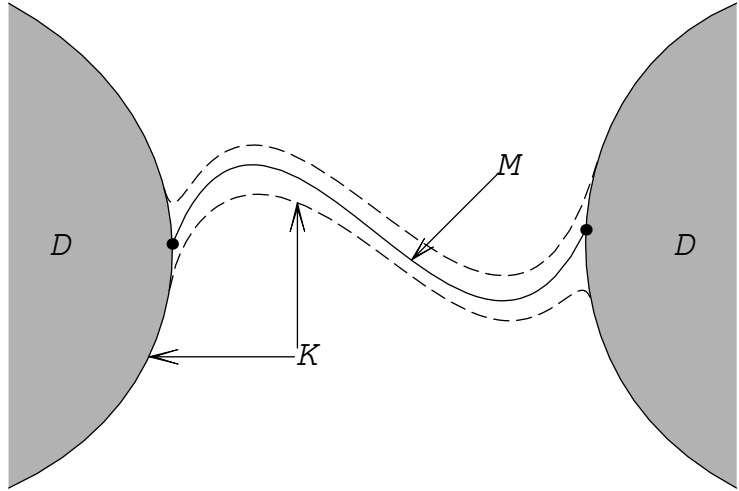}
  \caption{A handlebody $K$ with center $E = D \cup M$}\label{fig:HandleBody}
\end{figure}
\end{center}

Such $K$ will be called a {\it \spsc\ handlebody with center $E$}
(Figure \ref{fig:HandleBody}). The existence of a strong deformation 
retraction of $K$ onto $E$ implies that {\it $K$ is homotopically equivalent to $E$}. 

It is well known 
that any totally real submanifold $M$ in $\C^n$ (or in any complex manifold) 
has a basis of \spsc\ tubular neighborhoods. (If $M \subset\C^n$ is compact 
and of class $\cC^2$, we may take neighborhoods defined by the Euclidean distance 
to $M$.) Hence  the above problem is nontrivial only along the attaching 
submanifold $S=D\cap M \subset bD$. If $S$ fails to be Legendrian in $bD$ 
at some point $p\in S$ then $D\cup M$ may have a nontrivial local envelope 
of holomorphy at $p$, containing small analytic discs with boundaries in 
$bD\cup M$ (this follows from the results in \cite{Henkin}),
and in such case there exist no small pseudoconvex neighborhoods.
Local envelope may also appear at points $p\in M$ for which $T_p M$ 
contains a nontrivial complex subspace; see \cite{Bishop}. This justifies
the above hypotheses on $M$ and $S$.

The simpler problem concerning the existence of a basis of (strongly) 
pseudoconvex neighborhoods of $E=D\cup M$, without insisting on the
existence of a deformation retraction onto $E$, has been considered by several
authors; see e.g.\ Stolzenberg \cite{Stolzenberg}, H\"ormander and Wermer 
\cite{HormanderWermer}, Forn\ae ss and Stout \cite{FornaessStout1}, 
\cite{FornaessStout2}, Chirka and Smirnov \cite{SmirnovChirka}, 
and Rosay \cite{Rosay}. However, in many problems one actually needs 
\spsc\ handlebodies which have `the same shape' as $D\cup M$. 

An important general construction of handlebodies was given in 1990 by Eliashberg 
(Lemma 3.4.3.\ in \cite{Eliashberg}). 
Write the coordinates on $\C^n$ in 
the form $z=x+iy$, with $x,y\in\R^n$. Set $|x|^2=x_1^2+\ldots + x_n^2$, 
$|y|^2=y_1^2+\ldots+ y_n^2$. Let 
$$ 
	D_\l =\{x+iy \in \C^n\colon |y|^2 \ge 1+ \l |x|^2\},\quad 
	M = \{iy \colon |y|\le 1\}.
$$
Thus $M$ is the unit ball in the Lagrangian subspace $i\R^n\subset \C^n$,
attached to the quadric domain $D_\l$ along the $(n-1)$-sphere 
$S=bM=\{iy\colon |y|=1\} \subset bD_\l$ which is Legendrian in $bD_\l$. 
Note that $D_\l$ is \spsc\ precisely when $\l > 1$. 
In this situation, Lemma 3.4.3.\ in \cite{Eliashberg} gives for each 
open set $U\supset M$ a strongly pseudoconvex handlebody 
$K=\{|y|\ge \varphi(|x|)\}$, for a suitably chosen function $\varphi$,
such that $K \subset D_\l \cup U$ and the center of $K$ 
equals $E_\l=D_\l \cup M$. In \cite{Eliashberg} this  was used 
in the construction of Stein manifolds with prescribed homotopy type 
(see also Gompf \cite{Gompf} and Chapter 11 in \cite{GompfStipsicz}).

Five years later, in 1995, B.\ Boonstra \cite{Boonstra} (Ph.\ D.\ dissertation,
unpublished) constructed handlebodies whose center is the union of 
an ellipsoid with a Lagrangian plane in $\C^n$. He also constructed handlebodies in 
more general situations by the `osculation and patching' technique.  
Even though Boonstra cited Eliashberg's paper \cite{Eliashberg},
his construction seems independent since the details are somewhat different.

The content of this paper is as follows.
In section 2 we obtain a differential condition on a function $f$ 
which gives the necessary and sufficient condition for (strong) pseudoconvexity 
of the domain $D_+=\{x+iy \in \C^n \colon |y| \ge f(|x|)\}$ 
(resp.\ of $D_-=\{|y| \le f(|x|)\}$) along the hypersurface 
$\Sigma=\{|y|=f(|x|)\}$ (Proposition 2.1 and Corollary 2.2). 
A sufficient condition for strong pseudoconvexity of such domains was 
obtained earlier by Eliashberg; see (*, **) on p.\ 39 of \cite{Eliashberg}.
Our derivation of these conditions is different from the one in 
\cite{Eliashberg} and is somewhat similar to the one in 
\cite{Boonstra}. 

In section 3 we prove Proposition 3.1 which is the same as 
Lemma 3.4.3.\ in \cite{Eliashberg}. Our proof, based
on the differential conditions from Section 2, is  
similar to the original proof in \cite{Eliashberg}, but differs
from it in certain details. The extension to handles of lower 
dimension is immediate; see Lemma 3.1.1.\ in \cite{Eliashberg}.

Proposition 3.3 in the same section gives an explicit construction 
of \spsc\ handlebodies whose center is the union of
$D=\{x+iy\in \C^n\colon |y|^2 \le \l |x|^2+ 1\}$ (with $\l<1$)
and the Lagrangian plane $i\R^n$. 
Note that $D$ is \spsc\ precisely when $\l<1$; 
it is an unbounded hyperboloid when $0< \l<1$, a tube when $\l=0$, 
and a bounded ellipsoid when $\l<0$. Boonstra \cite{Boonstra} found explicit 
handlebodies for $\l< 0$ and gave an indirect construction for $0\le \l<1$. 
We give an explicit construction for all values $\l<1$. 

In Sect.\ 4 we construct monotone families of \spsc\ handlebodies whose 
center is the union of a sublevel set of a general quadratic \spsh\ function 
$\rho\colon \C^n\to\R$ and an attached disc $M \subset \R^k \subset \C^n$ 
that passes through the critical point $0\in\C^n$ of $\rho$. 
Unlike the handlebodies constructed by Eliashberg (or in section 3 above),
these handlebodies are not `thin' everywhere around $M$, but only in a 
smaller \nbd\ of the origin. Indeed these handlebodies are sublevel
sets of a certain noncritical \spsh\ function. 
The construction is independent from the one in \cite{Eliashberg} 
(and from the rest of this paper) and is much simpler. A crucial use of 
this result was made in \cite{Forstneric} (Lemma 6.7) in the
construction of \holo\ submersions of Stein manifolds to 
complex Euclidean spaces. 

Using standard bumping and patching techniques for strongly
pluri\-sub\-harmonic functions one may adapt the construction of handlebodies 
in \cite{Eliashberg} (and in this paper) to more general handle attachments,
assuming of course that the boundary of the handle is Legendrian in the boundary of 
the domain $D$. A particularly simple case is when the handle $M$ is 
real analytic along $bM$; in such case $M$ can be locally flattened
near any point $p\in bM$ by a local biholomorphic change of coordinates,
and the resulting domain can be osculated along $bM\subset D$ 
by a quadratic model of the type considered in \cite{Eliashberg} 
or in this paper. (This was used for instance in \cite{Rosay},
but with the weaker conclusion that $D\cup M$ admits a Stein neighborhood
basis. Certain cases have been treated by Boonstra \cite{Boonstra}, 
but his work remains unpublished.) 

The case of smooth (but non real-analytic) handles can possibly 
be handled by using coordinate changes near points 
$p\in bM$ which are $\overline\partial$-flat on $M$. Such coordinate 
changes clearly preserve strong pseudoconvexity of $bD$ locally
near $p$. However, to see that the model handlebodies remain strongly 
pseudoconvex under such coordinate changes, one must estimate the terms 
in their Levi form coming from the non-holomorphic terms in the 
coordinate change. We are not aware of any published work in this direction.

Professor Eliashberg informed us in a private communication 
(May 14, 2003) that a solution of the handlebody problem for 
handles $M$ of different topological type (i.e., non-disc type) 
follows from the case of disc-handles. Indeed, taking any 
Morse function on $M$ which is constant on $bM$, one decomposes $M$ 
into a union of disc-handles and then successively applies the 
disc-handle lemma. (The details do not seem to exist in print.)
We wish to thank him for this remark.

%
%
\section{Pseudoconvexity of spherical domains}\label{Pseudoconvexity}
Let $z=(z_1,\ldots,z_n) = x+iy \in\C^n$, with $z_k=x_k+iy_k$ for $k=1,\ldots,n$.
Set $|x|^2=x_1^2+\ldots + x_n^2$, $|y|^2=y_1^2+\ldots+ y_n^2$.
Let $U$ be a nonempty open set in $\R^n$ which is invariant under the
action of the orthogonal group $O(n)$ (i.e., $x\in U$ and $|x'|=|x|$ 
implies $x'\in U$).  Set $I =\{|x|^2 \colon x\in U \} \subset \R_+$. 
Assume that $\theta\colon I\to (0,+\infty)$ is a positive function of class 
$\cC^2$.

\proclaim 2.1 Proposition: Let $n>1$. The domain 
\begin{equation}\label{lab1}
	D_- = \{x+iy \in\C^n \colon x\in U,\ |y|^2 < \theta(|x|^2)\}   
\end{equation}
is strongly pseudo\-convex along the hypersurface $\Sigma=\{|y|^2 = \theta(|x|^2)\}$ 
\iff\ $\theta$ satisfies the following differential inequalities on $I$:
\begin{equation}\label{lab2}
	\theta'<1, \qquad   
	2|x|^2 \theta \theta'' < (1-\theta') 
	\left( |x|^2 \theta'^2 + \theta \right).                     
\end{equation}
($\theta$ and its derivatives are calculated at $|x|^2$).
The domain 
\begin{equation}\label{lab3}
      D_+ = \{x+iy \in\C^n \colon x\in U,\ |y|^2 > \theta(|x|^2)\}   
\end{equation}
is \spsc\ along $\Sigma$ \iff\ the reverse inequalities hold in (\ref{lab2}). 
If $\theta$ solves the differential equation
\begin{equation}\label{lab4}
   2|x|^2 \theta \theta'' = 
   (1-\theta')\cdotp \left( |x|^2 \theta'^2 + \theta \right)
\end{equation}
then $\theta'<1$ implies that $D_-$ is weakly pseudoconvex along $\Sigma$
while $\theta'>1$ implies that $D_+$ is weakly pseudoconvex along $\Sigma$.

\demo Proof: \rm
Set $\rho(x+iy) = |y|^2 - \theta(|x|^2)$. 
A calculation gives for $1\le j\ne k\le n$
\begin{eqnarray*}
 -\rho_{z_k} &=& x_k\theta' + iy_k, \\
              -2\rho_{z_k {\bar z}_k} &= & 2x_k^2\theta'' + \theta' -1, \\
              -2\rho_{z_j \bar{z}_k}  &= & 2x_jx_k \theta'', 
\end{eqnarray*}
where $\theta$ and its derivatives are evaluated at $|x|^2$.
The calculation of the Levi form of $\Sigma=\{\rho=0\}$ can be simplified 
by observing that $\rho$ is invariant under the action of the real orthogonal 
group $O(n)$ on $\C^n$ by $A(x+iy)=Ax+iAy$ for $A\in O(n)$.
Fix a point $p=r+is\in\Sigma$ $(r,s\in\R^n)$. After an orthogonal rotation
we may assume that $p=(x_1+iy_1,i\wt y_2,\ldots,i\wt y_n)$, with $x_1=|r| \ge 0$. 
Applying another orthogonal map which restricts to the identity on 
$\C\times\{0\}^{n-1}$ we may further assume that $p=(x_1+iy_1,iy_2,0,\ldots,0)$, 
where $y_1^2+y_2^2=|s|^2=\theta(x_1^2)$. At this point we have
$$ 
	\rho_{z_1}(p) = -x_1\theta' - iy_1, \quad
        \rho_{z_2}(p) = -iy_2, \quad 
        \rho_{z_k}(p) =0 \ {\rm for}\ k=3,\ldots,n. 
$$
Hence the complex tangent space 
$T_p^\C \Sigma = \{v\in \C^n \colon \sum_{k=1}^n {\di\rho\over \di z_k}(p)v_k=0\}$
consists of all $v\in\C^n$  satisfying 
$v_1=-\l iy_2$, $v_2=\l(x_1\theta'+iy_1)$ for arbitrary choices of 
$\l\in\C$ and $v''=(v_3,\ldots,v_n) \in\C^{n-2}$. We also have
\begin{eqnarray*}
	2\rho_{z_1\bar{z}_1}(p) &= & 1-\theta' -2x_1^2\theta'', \\  
        2\rho_{z_k\bar{z}_k}(p) &= & 1-\theta'  \qquad \qquad (k=2,\ldots,n), \\
        2\rho_{z_j\bar{z}_k}(p) &= & 0\quad\quad\qquad \qquad (1\le j\ne k\le n).
\end{eqnarray*}
For $v \in T^\C_p\Sigma$ we thus get (noting that $y_1^2+y_2^2 =\theta(x_1^2)$)
\begin{eqnarray}
 	2\cL_\rho(p;v) &=& \bigl(1-\theta' - 2x_1^2\theta''\bigr) |\l|^2 y_2^2
 	 + (1-\theta') |\l|^2 (x_1^2{\theta'}^2 + y_1^2) + (1-\theta')|v''|^2 \nonumber \\
 	 &= & |\l|^2 \left( -2x_1^2 y_2^2\theta'' + 
        (1-\theta')(x_1^2{\theta'}^2 + \theta) \right)  + (1-\theta') |v''|^2 \label{lab5}
\end{eqnarray}
where $\theta$ and its derivatives are evaluated at $x_1^2=|r|^2$.
Thus $\cL_\rho(p;v)>0$ for all choices of $\l\in\C$ and $v''\in\C^{n-2}$ 
with $|\l|^2+|v''|^2>0$ \iff\ 
$$
	\theta'<1, \qquad 
	2x_1^2 y_2^2\theta'' < (1-\theta')(x_1^2{\theta'}^2 + \theta).
$$
Observe that $0\le y_2^2\le |s|^2=\theta(x_1^2)$, and $y_2^2$ assumes 
both extreme values $0$ and $\theta(x_1^2)$ when $(y_1,y_2)$ traces the circle 
$y_1^2+y_2^2=\theta(x_1^2)$. Thus the second inequality above holds at all points 
of this circle precisely when it holds at the point $y_1=0$, $y_2=\sqrt{\theta(x_1^2)}$.
This gives the conditions
$$  
      \theta'<1,\qquad	2x_1^2 \theta\theta'' < (1-\theta')(x_1^2{\theta'}^2 + \theta) 
$$
characterizing strong pseudoconvexity of $D_-$ along the mentioned circle
in $\Sigma$. Since $x_1=|r|$, the above is equivalent to the pair of inequalities (\ref{lab2}) 
at $p=r+is$. Similarly we see that negativity of $\cL_\rho(p;v)$ for all choices 
of $\l$ and $v''$ (which characterizes strong pseudoconvexity of $D_+$) is 
equivalent to the reverse inequalities in (\ref{lab2}). 

Assume now that $\theta$ satisfies (\ref{lab4}). As before we reduce to the case 
$p=x+iy=(x_1+iy_1,iy_2,0,\ldots,0)$. From (\ref{lab5}) we obtain
\begin{eqnarray*}
	2\cL_\rho(p;v) &= & |\l|^2 \bigl(  -2x_1^2 y_2^2\theta'' + 2x_1^2\theta \theta'' \bigr)
	                + (1-\theta') |v''|^2 \\ 
	             &= & 2|\l|^2 x_1^2 y_1^2 \theta'' + (1-\theta') |v''|^2. 
\end{eqnarray*}
(We used $\theta(x_1^2)-y_1^2 = y_2^2$.) From (\ref{lab4}) we see that $\theta''$ is of 
the same sign as $1-\theta'$. Thus $\theta'<1$ implies $\cL_\rho(p;v) \ge 0$, with 
equality precisely when $v''=0$ and $0=x_1 y_1=x\cdotp y$. In this case
$D_-=\{\rho<0\}$ is weakly \psc\ along $\Sigma=\{\rho=0\}$, strongly pseudoconvex
on $\{x+iy\in \Sigma\colon x\cdotp y\ne 0\}$, and has one zero eigenvalue of the 
Levi form at each point of $\{x+iy\in \Sigma\colon x\cdotp y = 0\}$.
When $\theta'>1$ the analogous conclusions hold for $D_+=\{\rho>0\}$.
If $\theta'=1$ holds identically then 
$\rho=|y|^2-|x|^2+c= - \Re \bigl( \sum_{j=1}^n z_j^2 \bigr )+c$
is pluriharmonic.
\endpr

The second inequality in (\ref{lab2}) simplifies further in the 
variables $|x|, |y|$:
\proclaim 2.2 Corollary: Let $U\subset\R^n \bs\{0\}$ be an $O(n)$-invariant
open set and $f\colon I\to (0,+\infty)$ a $\cC^2$ function on 
$I=\{|x|\colon x\in U\}$. The domain 
$D_-=\{x+iy \in\C^n \colon x\in U,\ |y| < f(|x|)\}$ is 
strongly pseudo\-convex along the hypersurface 
$\Sigma=\{|y|=f(|x|)\}$ \iff\ 
\begin{equation}\label{lab6}
    {ff'\over |x|}<1, \qquad
    f\cdotp \left( f'' + {f'^3\over |x|} \right) < 1
\end{equation}
for all $x\in U$, where $f$ and its derivatives are evaluated at $|x|$. 
The domain $D_+=\{x+iy\colon x\in U,\ |y| > f(|x|)\}$ is strongly pseudoconvex 
along $\Sigma$ when the reverse inequalities hold in (\ref{lab6}).
If $f$ satisfies the differential equation
\begin{equation}\label{lab7}
	f\cdotp \left( f'' + {f'^3\over |x|} \right) = 1          
\end{equation}
then $ff'/|x| <1$ implies that $D_-$ is weakly pseudoconvex along
$\Sigma$ while $ff'/|x| >1$ implies that $D_+$ is weakly 
pseudoconvex along $\Sigma$. 

\demo Proof: Set $t=|x|>0$ for $x\in U$. The functions $f$ and $\theta$ 
are related by $f(t)^2=\theta(t^2)$. Differentiation gives 
$f(t)f'(t) =t\theta'(t^2)$ whence $\theta'<1$ 
is equivalent to ${ff'/t}< 1$. Another differentiation 
of $f(t)f'(t)=t\theta'(t^2)$ gives 
$$
	ff'' + f'^2= 2t^2 \theta'' + \theta' = 2|x|^2 \theta'' +  ff'/t.
$$
Hence $2|x|^2 \theta''= ff'' + f'^2-ff'/t$. Multiplying by $\theta=f^2$
we obtain the first line in the following display. In the second line
we used $|x|^2 \theta'^2=(t\theta')^2=f^2f'^2$:
\begin{eqnarray*}
  2|x|^2\theta\theta'' & = &
    f^2 \bigl( ff'' + f'^2 - {ff'\over t} \bigr), \\
    (1-\theta') \left( |x|^2 \theta'^2 + \theta \right) &= &
   \bigl(1-{ff'\over t}\bigr) (f^2 f'^2 + f^2) = 
   f^2\bigl( 1-{ff'^3\over t} + f'^2 - {ff'\over t} \bigr). 
\end{eqnarray*}
Comparing the two sides, dividing by $f^2>0$ and cancelling the common 
terms $f'^2 - ff'/t$ we see that the second inequality in (\ref{lab2}) is equivalent 
to $f(f''+f'^3/t)<1$. Similarly one treats the other cases.
\endpr

\ni\bf 2.3 Remarks. \rm
(A) The differential inequalities (\ref{lab2}) and (\ref{lab6}) 
are invariant up to the sign 
with respect to taking the inverses. More precisely, assume 
$\theta'(|x_0|^2) \ne 0$ for $x_0\in U$ and denote by $\tau$ the local inverse 
of $\theta$. At points where $\theta'>0$ the inequalities (\ref{lab2}) transform into 
the reverse inequalities for $\tau$:
$$ 
	\tau'>1,  \qquad  2|y|^2 \tau \tau'' > 
        (1-\tau') \left( |y|^2 \tau'^2 + \tau \right).                 
$$ 
On the other hand, near points where $\theta'<0$ the inequalities 
(\ref{lab2}) transform into the same inequalities for $\tau=\theta^{-1}$. 
This can be explained geometrically as follows. If $\theta'(|x_0|^2)>0$ then 
for $x$ near $x_0$ we have $|y|^2<\theta(|x|^2)$ \iff\ $|x|^2>\tau(|y|^2)$, 
and strong pseudoconvexity of the latter region is equivalent to the above 
inequality for $\tau$ according to Proposition 2.1. If $\theta'(|x_0|^2)<0$ 
then for $x$ near $x_0$ we have $|y|^2<\theta(|x|^2)$ \iff\ $|x|^2 < \tau(|y|^2)$, 
and pseudoconvexity is now characterized by (\ref{lab2}). Similarly the equations
(\ref{lab4}) and (\ref{lab7}) are invariant with respect to taking the inverses.

\ni (B) If $f(t)$ ($t\in \R$) is a function of class $\cC^1$ and piecewise $\cC^2$, 
we adopt the convention that $f$ satisfies the second inequality in (\ref{lab6}) 
at a point of discontinuity $t_0$ of the second derivative $f''$ when both the left 
and the right limit of $f''$ at $t_0$ satisfies it. (At endpoints we consider only the 
one sided limit.) A similar convention is adopted for (\ref{lab2}).
\endpr

\ni\bf 2.4 Example. \rm   
We illustrate the above by looking at model domains 
defined by the quadratic function  
$$   
	\rho_\l(z)=\rho_\l(x+iy) = \l |x|^2 - |y|^2 \qquad (\l\in\R,\ z\in\C^n) 
$$
which will be used in the following section. Setting 
$g_{\l,a}(t)= +\sqrt{\l t^2 + a}$ we have 
$\{\rho_\l < -a\}=\{x+iy\in\C^n\colon |y| > g_{\l,a}(|x|)\}$.
From  $\pmatrix {{\di^2\rho_\l \over \di z_j\di {\bar z}_k}} = {(\l -1)\over 2}I$ 
we see that $\rho_\l$ is \spsh\ when $\l>1$, strongly plurisuperharmonic when 
$\l<1$, and $\rho_1(x+iy)=|x|^2-|y|^2=\Re\left(\sum_{j=1}^n z_j^2 \right)$ 
is pluriharmonic. It is easily verified directly that $g_{\l,a}$ satisfies (\ref{lab6}) on  
$\{t\in \R_+\colon \l t^2+a\ge 0\}$ if $\l<1$, and it satisfies the reverse 
inequalities in (\ref{lab6}) if $\l>1$. If $\l\ne 0$ then $g=g_{\l,a}$ satisfies 
the differential equation $g \left( g''+ {g'^3 \over \lambda t} \right) =\lambda$.

%
%
%
%
\section{Strongly pseudo\-con\-vex han\-dle\-bodies}
\label{SPH}

In this section we find functions $f \colon I\to (0, +\infty)$
on intervals $I\subset \R_+ =[0,+\infty)$ which satisfy one of 
the following pairs of differential inequalities:
\begin{eqnarray}
  f\cdotp \left( f'' + {f'^3\over t} \right) < 1 
    & \quad{\rm and} & \quad {ff' \over t}< 1,                    \label{lab8}\\
    f\cdotp \left( f'' + {f'^3\over t} \right) > 1 
    & \quad{\rm and} & \quad {ff'\over t} > 1.                    \label{lab9}
\end{eqnarray}
If $f$ is of class $\cC^1$ and piecewise $\cC^2$ then at a point of 
discontinuity of $f''$ it should be understood that $f$ satisfies the 
first inequality in (\ref{lab8}) resp.\ in (\ref{lab9}) if the one-sided limits of $f''$ 
at that point satisfy it. By Corollary 2.2 the condition (\ref{lab8}) characterizes 
strong pseudoconvexity of the domain $\{x+iy\in\C^n\colon |y|<f(|x|)\}$ 
along $\Sigma=\{|y|=f(|x|)\}$ while (\ref{lab9}) does the same for
$\{x+iy \in \C^n \colon |y|>f(|x|)\}$.

\medskip \ni\bf 3.1 Proposition. \sl 
Let $\l>1, a>0$ and $g(t)=+\sqrt{\l t^2 + a}$. For every sufficiently 
small $\e>0$ there exists a number $\s=\sigma(\e)\in (0,\e)$ and a 
continuous, positive, strictly increasing function 
$f=f_\e \colon [\s,+\infty) \to [f(\s),+\infty)$ which is 
$\cC^\infty$ on $(\s,+\infty)$, satisfies (\ref{lab9}) and also the following:
\begin{itemize}
\item[(i)]   $f(t)= g(t)$ for $t\ge \e$, 
\item[(ii)]  $f(t)< g(t)$ for $\sigma\le t<\e$,  
\item[(iii)] $f'(\s+) = \lim_{t\downarrow \s} f'(t)=+\infty$, and
\item[(iv)]  the inverse function $f^{-1} \colon \R_+\to [\s,+\infty)$ is of 
class $\cC^\infty$ and satisfies (\ref{lab8}) provided that we set 
$f^{-1}(u)= \sigma$ for $0\le u\le f(\sigma)$.
\end{itemize}
\medskip\rm

\begin{center}
\begin{figure}[hbt]
  \centering\epsfig{file=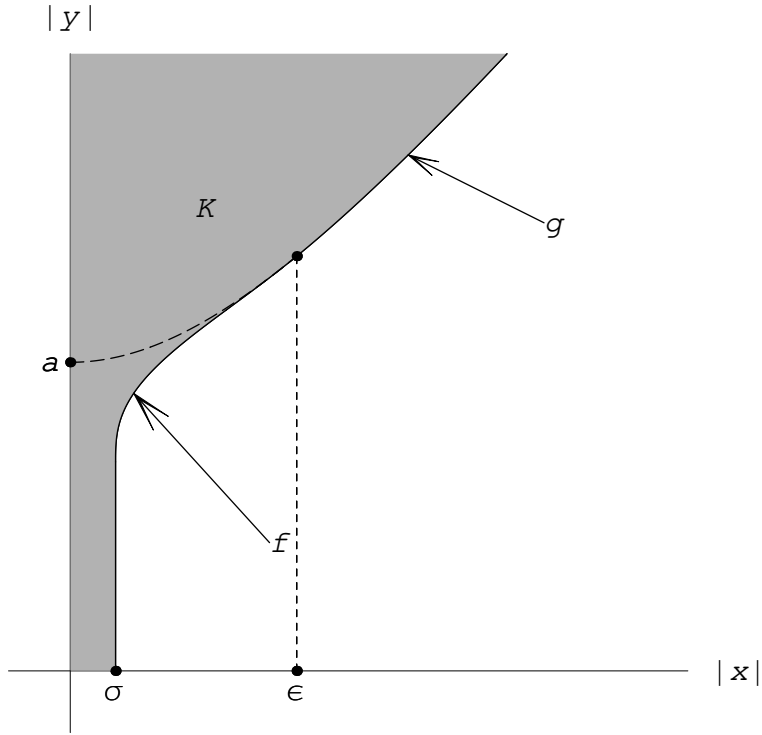}
  \caption{The handlebody $K$}\label{fig:Ke}
\end{figure}
\end{center}

\proclaim 3.2 Corollary:
Let $\l>1, a>0$, $D=\{x+iy\in\C^n \colon |y|^2 \ge \l |x|^2 + a\}$
and $M=\{iy \colon y\in\R^n,\ |y|\le a\}$. If $f$ satisfies 
Proposition 3.1 then $K = \{x+iy\in\C^n \colon |x| \le f^{-1}(|y|)\}$ 
(Figure \ref{fig:Ke}) is a smooth \spsc\ handlebody with center $E=D\cup M$,
satisfying $D\cup \{|x| \le  \s\} \subset K \subset  D\cup \{|x| < \e\}$.

%
%

\ni\it Remark. \rm  We have already said in the Introduction
that Proposition 3.1 (and Corollary 3.2) is the same as 
Lemma 3.4.3.\ in \cite{Eliashberg}. The handlebodies 
on figures  \ref{fig:Ke} and \ref{fig:L} are shown in the 
coordinate system $(|x|,|y|)\in \R_+^2$; the actual handlebody 
is the preimage under the map $x+iy \to (|x|,|y|)$. 

\begin{center}
\begin{figure}[hbt]
  \centering\epsfig{file=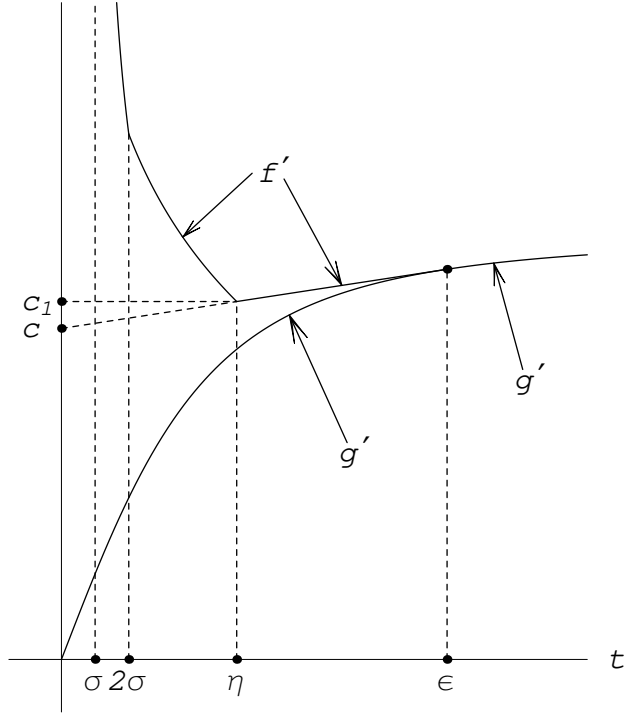}
  \caption{The graph of $f'$}\label{fig:fcrtaK}
\end{figure}
\end{center}

\begin{center}
\begin{figure}[hbt]
  \centering\epsfig{file=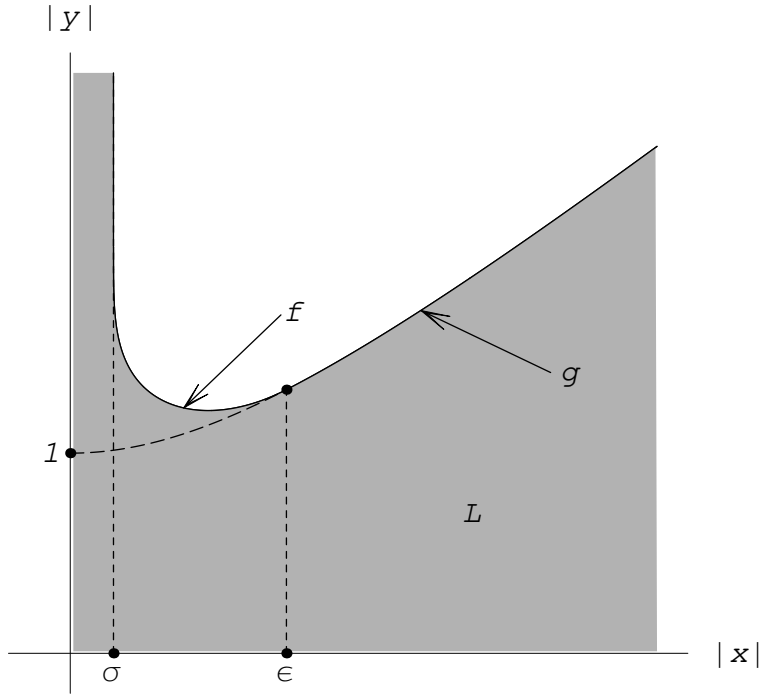}
  \caption{The handlebody $L$}\label{fig:L}
\end{figure}
\end{center}

\begin{center}
\begin{figure}[hbt]
  \centering\epsfig{file=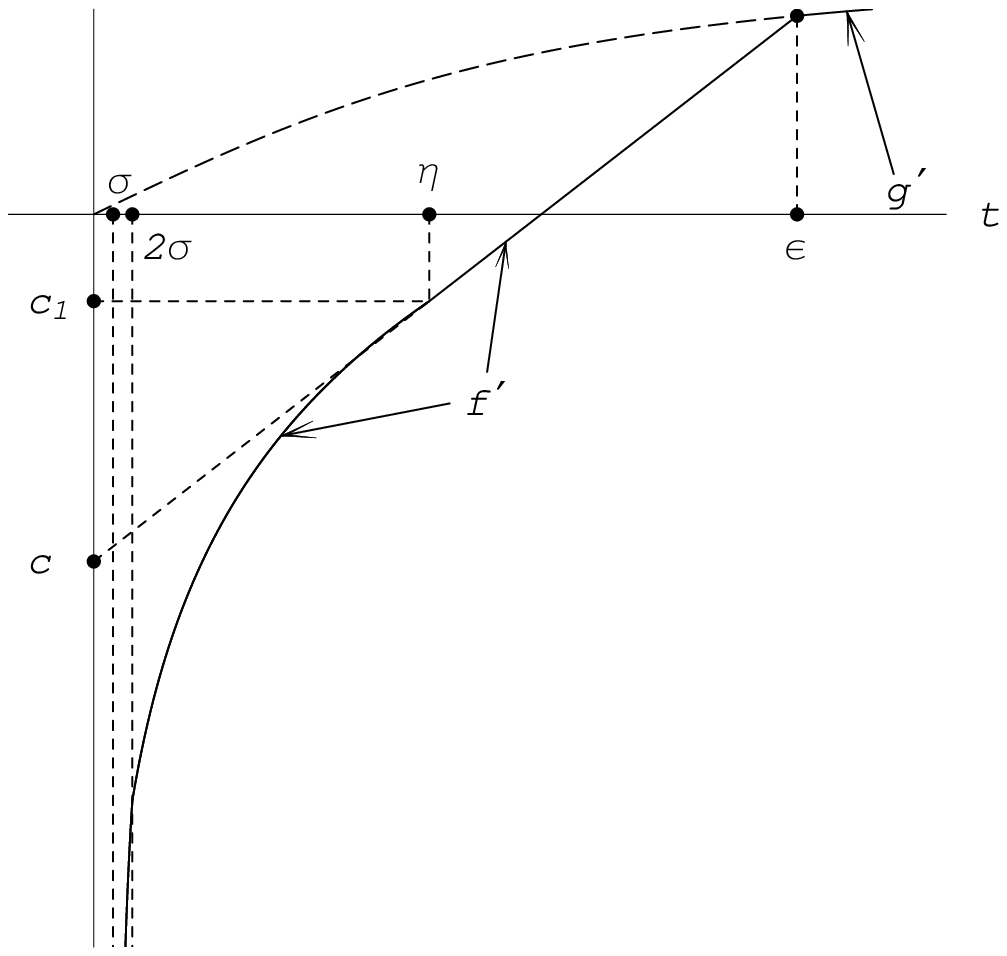}
  \caption{The graph of $f'$}\label{fig:fcrtaL}
\end{figure}
\end{center}

\demo Proof of Proposition 3.1: 
Without loss of generality we may take $a=1$ and 
$g(t)=\sqrt{\l t^2+1}$ (the general case follows by rescaling). 
A calculation gives for $t>0$
$$ 
	g'(t)  = {\l t\over\ g(t)}>0,\quad 
	g''(t) = {\l\over g(t)^3}>0, \quad 
        g'''(t)= -{3\l^2 t\over g(t)^5}<0   
$$
which shows that $g$ is increasing, convex, and $g'$ is concave. 
We also obtain $g'(t)-tg''(t)=\l^2 t^3/g(t)^3>0$.
Fix a small $\e>0$ and let $c:=g'(\e)-\e g''(\e)>0$. Choose a number 
$0<\eta< \min (\e,c^3/3)$ and let $c_1:=c+\eta g''(\e)$. Let $\s>0$ be a 
number satisfying $2\s<\eta<\e$ (its precise value will be determined later). 
We shall first obtain a solution $f$ of class $\cC^1$ and piecewise $\cC^2$ 
on $(\sigma,+\infty)$; the final solution will be obtained by smoothing. Let  
$$ 
	f(t)= g(\e)+ \int_{\e}^t f'(\tau)\, d\tau  \qquad (\s \le t< +\infty)
$$
where $f'$ is a continuous and piecewise $\cC^1$ function defined as follows: 
$$              
             f'(t)=\cases{ g'(t),                &if $\e \le t$;\cr
                          g'(\e)+g''(\e)(t-\e), &if $\eta \le t < \e$;          \cr
                          c_1 + \eta\log(\eta/t), &if $2\sigma \le t < \eta$;     \cr
                          {2\sqrt\sigma / \sqrt{t-\sigma}}, &if $\s<t< 2\s$.\cr
                          }
$$
The graph of $f'$ is shown on Figure \ref{fig:fcrtaK}.
(However, due to technical difficulties we show the case for 
large $\e$. For small $\e>0$ the derivative of the linear part 
of the graph should be close to $\l >1$. The same remark applies to 
Figure \ref{fig:fcrtaL}.) 

Note that $f'$ is continuous at $t=\eta$, with $f'(\eta)=c_1$. To insure the 
continuity of $f'$ at $t=2\s$ we choose $\s$ to be the solution 
of $c_1 + \eta\log(\eta/2\s)=2$.
Clearly $f(t)=g(t)$ for $t\ge \e$. It is also clear that $f'(t) >g'(t)$ 
for $\s< t< \e$: on $t\in [\eta,\e]$ the graph of $f'$ is the tangent line 
to the graph of $g'$ at $(\e,g'(\e))$ which stays above $g'$ due the
to concavity of $g'$; on $(\s,\eta]$ this is clear since 
$g'$ is increasing while $f'$ is decreasing. Hence $f$ is strictly 
increasing and satisfies $f(t)<g(t)$ for $\s \le t<\e$. Also 
$f'(\s+)=+\infty$. It remains to show that $f(\s+)>0$ and that $f$ 
satisfies (\ref{lab9}) on $(\s,\e)$. 

%
%
\ni \bf Case 1: \rm $\eta\le t < \e$. On this interval 
$$ 
	f'(t)= g'(\e)+g''(\e)(t-\e) = c+ tg''(\e) > tg''(\e).
$$
The graph of $f'$ is the tangent line to the graph of $g'$ at 
the point $(\e,g'(\e))$. Since $g'$ is strongly concave, we conclude 
$f'(t) > g'(t)$ for all $t\in [\eta, \e)$. We have
$$ 
	f(t) >  g(\e) - \int_0^\e \bigl(g'(\e)+g''(\e)(\tau -\e)\bigr) d\tau
	     >  g(\e) - \e g'(\e) = 1/g(\e).
$$
Since $f'(t)>0$ and $f''(t)=g''(\e)$, we get
$f\left( f''+f'^3/t \right)  > ff'' > {g''(\e)/g(\e)} = \l/g(\e)^4$
which is $>1$ if $\e$ is small (since $\l>1$ and 
$g(\e)\approx g(0)=1$). From $f(t) > 1/g(\e)$ 
and $f'(t) > tg''(\e)$ we also get $f(t)f'(t)/t > g''(\e)/g(\e) >1$. 
\hfill\break\smallskip\break
\ni \bf Case 2: \rm $2\s \le t < \eta$. Using $f(\eta)>1/g(\e)$, 
$f'(\eta)< g'(\e)$ (Case 1) we get
$$ 
    f(t) > f(\eta) - \int_0^\eta \left( f'(\eta) -\eta\log(\tau/\eta)\right)d\tau 
                  >   {1 \over g(\e)} - \eta g'(\e)-\eta^2 =:M. 
$$
Clearly $M>1/2$ when $\e>0$ is small. From $f''(t)=-\eta/t$, 
$f'(t)=f'(\eta)+\eta \log(\eta/t)> f'(\eta) > c>0$
and $0<3\eta<c^3$ we obtain
$$  
	f\left(f''+ {f'^3 \over t} \right)- 1 > 
         M\left( {-\eta \over t} + {c^3\over t} \right) -1
         > {M\over t}\left( c^3-\eta - 2t \right)
         > {M\over \eta}\left( c^3- 3\eta \right)>0.
$$
Also, $ff'/t > Mf'(\eta)/\eta>Mc/\eta> 3M/c^2>1$ (since $c>0$ is small)
which verifies the second inequality in (\ref{lab9}). 
\hfill\break\smallskip\break
\ni \bf Case 3: \rm $\s < t < 2\s$. As before we easily obtain a lower bound
$f(t)>1/2$ provided that $\e>0$ is sufficiently small. We have
$f'(t)={2\sqrt\sigma / \sqrt{t-\sigma}}$, 
$f''(t)=-\sqrt\s / \sqrt{t-\s}^3$, and hence
$$  
	f\left(f''+ f'^3/t \right)  > {1\over 2} 
                 \left( {-\sqrt\s \over \sqrt{t-\s}^3}  + 
                 {8\s \sqrt\s \over t \sqrt{t-\s}^3} \right) 
        > {\sqrt\s \over 2\sqrt{t-\s}^3 } (-1+4) \ge {3\over 2\s} >1. 
$$
The second inequality in (\ref{lab2}) is trivial as in Case 2. 
\hfill\break

The function $f$ constructed above is invertible and its inverse 
function $f^{-1}\colon [f(\s),+\infty) \to [\s,+\infty)$ 
is of class $\cC^1$, piecewise $\cC^2$ (actually piecewise 
\ra), and satisfies (\ref{lab8}). We extend $f^{-1}$ to $[0,+\infty)$ 
by  taking $f^{-1}(u)=\s$ for $u\in [0,f(\s)]$; this extension 
satisfies the same properties also near the point $u=f(\s)$. The final 
solution is obtained by smoothing $h:=f^{-1}$ in a small \nbd\ of any 
point of discontinuity of it second derivative. (We interpolate smoothly 
between the left and the right limit of $h''$ at such a point and 
integrate twice to obtain the new $h$. This does not change $h$ and 
$h'$ very much and hence the inequality (\ref{lab8}) is preserved.) 
This completes the proof. 
\endpr

%
%
%
%
A small modification of the above construction gives \spsc\ handlebodies 
$L\subset \C^n$ with center 
$$  
	E = \{x+iy\in\C^n \colon |y|^2 \le \l |x|^2+ 1\} \cup i\R^n
        \qquad (\l<1). 
$$
%
%
A typical $L$ is shown on Figure \ref{fig:L}.
Observe that $D=\{|y|^2 \le \l |x|^2+ 1\}$ is \spsc\ precisely when $\l<1$. 
It is an unbounded hyperboloid when $0< \l<1$, a tube when $\l=0$ 
and a bounded ellipsoid when $\l<0$. The Lagrangian plane $i\R^n$
is an (unbouded) handle attached to $D$ along the sphere 
$\{iy \colon y\in \R^n,\ |y|=1\}$. Boonstra \cite{Boonstra} found explicit 
handlebodies for $\l< 0$ and gave an indirect `bumping and patching' 
construction for $0\le \l<1$. We give an explicit construction 
for all $\l<1$. (Our example is easily modified to obtain 
handlebodies with center $D\cup M$ where $M \subset i\R^n$ is a 
compact domain such that $D\cap i\R^n$ is contained in the 
relative interior of $M$.) Set  
$$
	L= \{x+iy\colon |x| > \s,\ |y| \le f(|x|\} 
	\cup \{x+iy\colon |x|\le \s\}
$$
where $f$ is given by the following proposition.
\endpr

\medskip \ni\bf 3.3 Proposition. \sl 
Let $\l<1$ and $g(t)=+\sqrt{\l t^2 + 1}$. For every sufficiently 
small $\e>0$ there exists a number $\s=\sigma(\e)\in (0,\e)$ and a 
continuous function $f\colon [\s,+\infty) \to (0,+\infty)$, 
smooth on $(\s,+\infty)$, which satisfies the inequalities 
(\ref{lab8}) and the following:
\begin{itemize}
\item[(i)]   $f(t)= g(t)$ for $t\ge \e$, 
\item[(ii)]  $f(t)> g(t)$ for $\sigma\le t<\e$,  
\item[(iii)] $f'(\s+) = \lim_{t\downarrow \s} f'(t)=-\infty$,  
\item[(iv)]  there exists a smooth inverse function $f^{-1}$ near 
the point $u=f(\s)$, with $f^{-1}(u)=\s$ for $u\ge f(\s)$, 
satisfying the inequalities (\ref{lab8}) on its domain. 
\end{itemize}
\medskip\rm

\demo Proof: 
Choose numbers $0<\s < 2\s< \eta<\e$; additional conditions will be imposed later. 
We have $g'(\e)/\e = \l/g(\e) <1$. Choose a number $k$ satisfying 
$g'(\e)/\e < k <1$ and let $c:=g'(\e)-k\e <0$. Clearly $c>-1$ if $\e$ is small.
Choosing $\eta>0$ sufficiently small we have 
$c_1:=g'(\e)+k(\eta-\e)= c+k\eta <0$ and $\eta+c_1^3<0$. 
Let $\s\in (0,\eta/2)$ solve $c_1-\eta\log(\eta/2\s)=-2$. With these choices we 
define $f$ on $(\s,+\infty)$ by $f(t)=g(\e)+ \int_{\e}^t f'(\tau)\, d\tau$ where 
$$              
             f'(t)=\cases{ g'(t),                &if $\e \le t$;\cr
                          g'(\e)+k(t-\e),        &if $\eta \le t < \e$;      \cr
                          c_1 - \eta\log(\eta/t), &if $2\sigma \le t < \eta$;     \cr
                          -{2\sqrt\sigma / \sqrt{t-\sigma}}, &if $\s<t< 2\s$.\cr
                          }
$$
%
%
The graph of $f'$ is shown on Figure \ref{fig:fcrtaL}.
We verify that $f$ satisfies (\ref{lab8}). For $t\ge \e$ this is 
clear since $f(t)=g(t)$. For $\eta\le t<\e$ we have
$$ 
	g(t)< f(t)\le g(\e)+ \big| \int_t^\e \left(g'(\e)+k(\tau -\e) \right) d\tau \big| 
	< g(\e) + \e(1 + |g'(\e)|).   
$$
By our choice of $k$ the graph of $f$ lies below the secant line through
$(0,0)$ and $(\e,g'(\e))$, and the secant is below $g'$ due to concavity
of $g'$. This gives $f'(t) < g'(\e)t/\e = \l t/g(\e)$. Also, $f''(t)=k$. 
At points $t\in [\eta,\e)$ where $f''(t)+f'(t)^3/t >0$ we thus have
$$ 
	f(t)\bigl(f''(t)+f'(t)^3/t \bigr) 
	\le \left( g(\e)+ \e(1 + |g'(\e)|) \right)  
	\left( k + \l^3 \e^2/g(\e)^3 \right) < 1
$$
provided that $\e>0$ is sufficiently small (since $k<1$, $g(\e)\approx 1$
and the other quantities are $O(\e)$). At points where 
$f''(t)+f'(t)^3/t \le 0$ the same estimate holds since $f(t)>0$. Also, 
$f(t)f'(t)/t < \bigl(g(\e)+O(\e)\bigr) \l/g(\e) =\l +O(\e) < 1$ 
if $\e>0$ is small. 

For $t\in [2\s,\eta]$ the estimates (\ref{lab8}) are almost trivial:
from $f'(t)\le f'(\eta)=c_1<0$ and $f''(t)=\eta/t$ we get
$f''(t)+f'(t)^3/t \le (\eta +c_1^3)/t < 0$ which implies the first 
estimate in (\ref{lab8}) (since $f(t)>g(t)>0$). Also $f(t)f'(t)/t < 0$ and hence
(\ref{lab8}) holds. Similarly we verify (\ref{lab8}) on $(\s,2\s]$. We complete
the proof as in Proposition 3.1 by smoothing $f^{-1}$. 
\endpr

%
%
%
%
\section{Handlebodies on general quadratic domains}
\label{HandlebodiesGQD}

In this sec\-tion we con\-sider han\-dle\-bod\-ies mod\-eled on gen\-eral quadratic
\spsh\ functions $\rho\colon\C^n\to\R$. Choose a $k\in\{0,1,\ldots, n\}$
and write the coordinates on $\C^n$ in the form $\z=(z,w)$, with $z=x+iy \in\C^k$ 
and $w=u+iv\in\C^{n-k}$. Let $A, B$ be positive definite real symmetric 
matrices of  dimension $k$ resp.\ $n-k$. Denote by $\langle \cdotp,\cdotp\rangle$
the Euclidean inner product on any $\R^m$. Given these choices let
\begin{equation} \label{lab12}
	 \rho(z,w) = Q(y,w)-|x|^2,
	 \quad Q(y,w)=\langle Ay,y\rangle + \langle Bv,v\rangle + |u|^2.    			     						
\end{equation}
It is easily seen that $\rho$ is \spsh\ \iff\ all eigenvalues of $A$ are larger than
$1$. (Equivalently, the matrix $A-I$ must be positive definite which we denote
by $A>I$.)  Clearly $\rho$ has a Morse critical point of index $k$ at the origin 
and no other critical points. It is proved in \cite{HarveyWells} that every 
Morse critical point of a \spsh\ function is of this form in some local 
holomorphic coordinates, modulo terms of order $>2$.
 
Assume now that $k\ge 1$. Let $\Lambda^k = \{(x+i0,0)\in\C^n \colon x\in \R^k\}$.
We identify $x\in\R^k$ with $(x+i0,0)\in \Lambda^k\subset \C^n$ 
when appropriate. 

\medskip\ni\bf 4.1 Proposition. \sl {\rm (Notation  as above.)}
Let $\rho$ be given by (\ref{lab12}) where $A>I$, $B>0$. 
Given $r>0$, $\e>0$ there exist constants $0<r<c_0<R$, $\d>0$ and a smooth, 
increasing, weakly convex function $h\colon \R_+\to\R_+$ 
such that $\tau(z,w)\colon=Q(y,w)-h(|x|^2)$ is a \spsh\ function 
on $\C^n$, with a Morse critical point of index $k$ 
at $0\in\C^n$, satisfying 
\begin{itemize}
\item[(i)]   for $|x|^2\le r$ we have $\tau(z,w)=Q(y,w)-\d |x|^2$,
\item[(ii)] for $|x|^2\ge R$  we have $\tau(z,w)=\rho(z,w)+c_0$, and
\item[(iii)] $\{\rho \le -c_0\}\cup \Lambda^k \subset \{\tau\le 0\} 
\subset \{\rho<-r\}\cup \{Q<\e\}$.
\end{itemize}

%
%

\proclaim 4.2 Corollary: 
For every sufficiently small $c>0$ the  set 
$K_c =\{\tau \le c\}$ is a \spsc\ handlebody with center
$$
	E_{c-c_0} =\{\z\in\C^n\colon \rho(\z)\le c-c_0 \} \cup \Lambda^k,
$$
satisfying $E_{c-c_0}\subset K_c\subset  \{\rho<-r\}\cup \{Q<\e\}$
(Figure 6).

\demo Proof of Proposition 4.1: 
We modify slightly the construction in Lemma 6.7 of \cite{Forstneric} 
(the function constructed there was not Morse). 
Let $t_0=r+\e$. Choose $0<\d <1, \mu>1$ such that 
$\d t_0<\e$ and $1<\mu+\d<\l_1$ where $\l_1>1$ denotes the smallest eigenvalue
of $A$. Set $R= {\mu^2 t_0 /(\mu+\d-1)^2}$ and 
$$      h(t) = \cases{ \d t,& if $0\le t\le t_0$; \cr
                       \d t + \mu(\sqrt t -\sqrt t_0)^2,& if $t_0<t\le R$; \cr
                       t -R +h(R),& if $R<t$. \cr}
$$

\begin{center}
\begin{figure}[hbt]
  \centering\epsfig{file=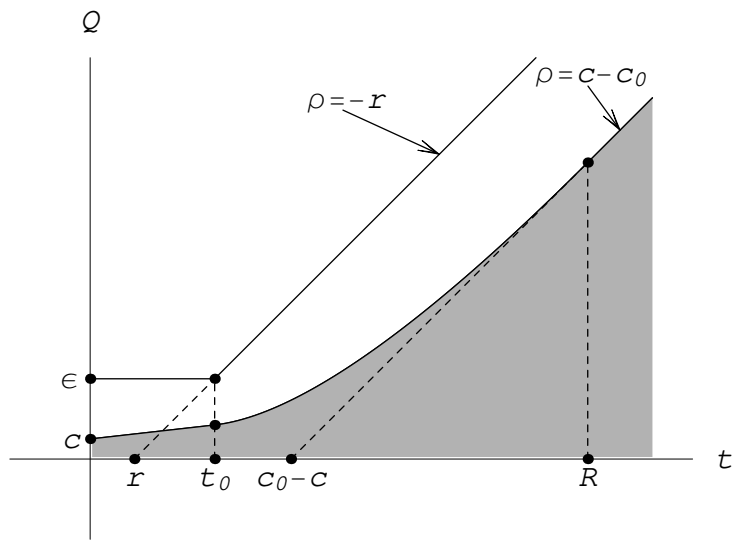}
  \caption{The handlebody $K_c = \{\tau \le c \}$}\label{fig:Kc}
\end{figure}
\end{center}

It is easily verified that $h$ is an increasing convex function of class $\cC^1$ 
and piecewise $\cC^2$ on $\R$ which satisfies 
$$ 
	2t\ddot h + \dot h=\mu+\d <\l_1 \qquad (t_0\le t\le R) 
$$
and $\d = \dot h(t_0) \le \dot h(t) \le 1=\dot h(R)$ for all $t\in\R$.
By smoothing $h$ we obtain an increasing convex $\cC^\infty$ function, still denoted $h$, 
which equals $\d t$ for $0\le t\le t_0$, it equals $t-R+h(R)$ for $t\ge R$, and
satisfies
$$ 
	\dot h(t)<\l_1,\quad  2t\ddot h(t) + \dot h(t) <\l_1 \qquad(t\in\R).
$$
A simple calculation shows that, as a consequence of these inequalities, the associated 
function $\tau$ is \spsh\ on $\C^n$ and satisfies Proposition 4.1 with $c_0=R-h(R)$.
(See the proof of Lemmas 6.7 and 6.8 in [F] for the details of this 
calculation.)
\endpr

%
%
%
%

\bigskip
\ni Institute of Mathematics, Physics and Mechanics,
University of Ljubljana, Jadranska 19, 1000 Ljubljana, Slovenia

\end{document}